\documentclass{svproc}
\usepackage{graphicx}%
\usepackage{multirow}%
\usepackage{amsmath,amssymb,amsfonts}%
\usepackage{mathrsfs}%
\usepackage[title]{appendix}%
\usepackage{xcolor}%
\usepackage{textcomp}%
\usepackage{manyfoot}%
\usepackage{booktabs}%
\usepackage{algorithm}%
\usepackage{algorithmicx}%
\usepackage{algpseudocode}%
\usepackage{listings}%
\usepackage{url}

\usepackage{epsfig}
\usepackage[hang]{caption2}
\usepackage[normalsize,it,hang]{subfigure}

\begin{document}
\mainmatter              
\title{{\tt \footnotesize The 6th International Symposium on Complex Systems (ISCS 2026) \newline June 03-05, 2026, La Rochelle, France} \vspace{1.4cm} \\
Sterile mosquito release \\ via intelligent proportional controllers}
\titlerunning{Sterile Insect Technique}  
%
\author{C\'{e}dric Join\inst{1} \and Luis Almeida\inst{2}
\and
Michel Fliess\inst{3,4}}
\authorrunning{C\'{e}dric Join et al.} 

\institute{Centre de Recherche en Automatique de Nancy, CNRS, Universit\'{e} de Lorraine, 54506 Vand{\oe}uvre-l\`{e}s-Nancy, France\\
\email{cedric.join@univ-lorraine.fr}
\and
{Sorbonne Universit\'{e}, Universit\'{e} Paris Cit\'{e}, CNRS, Laboratoire de Probabilités, Statistique et Modélisation (LPSM), 75005 Paris, France}\\
\email{luis.almeida@cnrs.fr}
\and
Laboratoire Jacques-Louis Lions, CNRS, Sorbonne Universit\'e, \\ 75005 Paris, France \\
{\tt michel.fliess@sorbonne-universite.fr, michel.fliess@swissknife.tech}
\and
AL.I.E.N., 7 rue Maurice Barr\`{e}s, 54330 V\'{e}zelise, France\\
{\tt michel.fliess@alien-sas.com}}

\maketitle

\begin{abstract}
The Sterile Insect Technique (SIT) against insect pests and insect vectors
consists of releasing males that have been previously sterilized in order to reduce or eliminate a specific
wild population. We study this complex control question via model-free control, ultra-local models, and intelligent proportional controllers that have already proven their effectiveness in various fields. They permit addressing, perhaps for the first time, the essential sampling question. Computer simulations are displayed and discussed. 
\end{abstract}

\newpage

\section{Introduction}\label{sec1}
%


Among the numerous ontological questions regarding the very nature of \textit{complex control systems}, let us single out a recent and enlightening review paper \cite{bo} in which it is suggested to relate it to artificial neural networks, deep learning, and generative artificial intelligence. His citation of a publication over $50$ years old by Ivakhnenko \cite{iva}, a famous Soviet Academician, reads: \newline \textit{Modern control theory, based on differential equations,
is not an adequate tool for solving the problems of
complex control systems. Constructing differential
equations to trace input-output paths requires a deductive, deterministic approach. However, this approach is
impractical for complex systems due to the difficulty
of identifying these paths}. \newline This communication, which follows those precepts, is based on \emph{model-free control} \cite{mfc1,mfc2}:
\begin{itemize}
    \item it is easy to implement \cite{nice};
     \item a simplistic \emph{ultra-local model} is combined with an \emph{intelligent proportional controller}.
\end{itemize}
 See \cite{ia} for preliminary remarks on connections with AI. Some thorough comparisons (see, e.g., \cite{madrid}) conclude with the superiority of this setting over other control viewpoints. Only a few references \cite{vila,diab,bara,kiz,precup,viet,ital,michel,exosk,alg,epilepsy,parkinson} in biomedicine and bioengineering are listed.

Mosquito-borne diseases (for instance, dengue, Zika, and chikungunya)  will be at the core of this paper: 
 \begin{itemize}
     \item It is the first time that a complex issue at the intersection of ecology and medicine will be investigated using the aforementioned tools.
     \item It is among the most important public health problems worldwide and is likely to continue rising in the near future, as many mosquito species expand their geographical ranges due to global trade and global warming. 
 \end{itemize}
Given past difficulties in vaccine and drug development, present mosquito-borne disease control programs rely on conventional mosquito reduction methods that are resource intensive and often not sustainable. In particular, insecticides pollute the environment, and their broad spectrum severely impacts biodiversity. Moreover, mosquito vectors rapidly develop insecticide resistance, which reduces control effectiveness. 
Consequently, there is an urgent need to develop robust, species-specific, and environmentally sustainable control strategies. Mechanical interventions, such as removing potential breeding sites and deploying mosquito traps, offer partial solutions but are usually not sufficient to control the mosquito population, and their implementation relies heavily on community adherence. This has led to the development of several biological control approaches that look very promising.

In this work, we will consider the {\em Sterile Insect Technique} ({\em SIT}), which consists of the release of large amounts of sterile male mosquitoes to reduce the
proliferation of mosquitoes by perturbing their reproduction. The technique has been successfully used in the field, not only for mosquitoes but also for other pests (see, for instance, \cite{ABBCD,dyck2021sterile}). Unlike the use of insecticides, this technique has the advantage of being species-specific and environmentally friendly, and is thus favored by many researchers and institutions.

In the existing literature, mathematical modeling plays an essential role in the safe and effective deployment of mosquito control strategies. Most of these models consist of systems of ordinary differential equations (ODEs) in order to capture the temporal dynamics of mosquito populations \cite{ABBCD,APSV}, but there also exist mathematical models that describe the dynamics of the invasive fronts of mosquito populations \cite{MEL}, as well as other models aimed at identifying ways to use SIT interventions to disrupt their propagation \cite{almeidaBarrier2022,almeidaRolling2023}.  
Emerging approaches in mosquito population control, such as the {\it rolling carpet} strategy \cite{almeidaRolling2023,almeidaRolling2025}, the use of optimal and feedback control \cite{almeida2022optimal,AAC,AAC2}, as well as some monotonicity techniques \cite{bliman} aim to improve the efficiency and robustness of interventions, keeping the system as close as necessary to the mosquito-free state. 

The various model-based approaches are quite intricate and, therefore, difficult to implement in practice. In this study, we intend: 
\begin{itemize}
    \item to validate, to some extent, i.e., without concrete experiments, the model-free control strategy obtained via a simple model-based description; 
    \item to check its robustness with respect to model mismatches, such as poor parameter knowledge. 
    \item to address the essential sampling question. In practice, sterile male mosquitoes are not released continuously (all the time) but only at certain moments: typically twice a week for Aedes mosquitoes (see \cite{gato2021sterile}). Therefore, it is natural to work in an impulsive control setting, for which our approach is particularly well suited.
\end{itemize}

The proposed control strategy is evaluated using several simulations. 
In the nominal case with pulses every 3 days, the reconstruction of the continuous command is satisfactory, and the reference trajectory is tracked perfectly. When we increase the duration between treatments to 6 days, the reference trajectory tracking degrades slightly because saturation is reached, but tracking remains very accurate overall. When we test robustness to sterile males dynamic uncertainty, performance remains very good even though the continuous-discrete representation no longer coincides perfectly. Finally, we perform 100 simulations with random parameter variations between 0.7 and 1.3 times the nominal values, and tracking is not guaranteed in only two cases where control is saturated at all times.

Sect. \ref{I} summarizes model-free control. A mathematical model is given  in Sect. \ref{descrip} in order to test our model-free control strategy via computer experiments. A realistic, i.e., discontinuous, setting is proposed in Sect. \ref{algo}. Sect. \ref{sec2} provides a quite detailed analysis of numerous numerical simulations. Some concluding remarks may be found in Sect. \ref{conclusion}.

\section{Model-free control}\label{I}

\subsection{Ultra-local model}\label{loop}
It has been demonstrated \cite{mfc1} (see \cite{bara} for a more comprehensible presentation) that under quite weak assumptions, any single-input single-output (SISO) system may be approximated by an \emph{ultra-local model}, which needs to be continuously updated,
\begin{align}\label{eq-1}
y^{(\nu)} = F + \alpha u,
\end{align}
where
\begin{itemize}
\item the control and output variables are respectively $u$ and $y$;
\item the derivation order $\nu$ may often be set to $1$; 
\item the constant $\alpha \in \mathbb{R}$, which is fixed in such a way that $\alpha u$ and
$y^{(\nu)}$ are of the same order of magnitude, does not need to be precisely computed;
\item $F$ subsumes not only the system dynamics but also external disturbances.
\end{itemize}

{Eq. \eqref{eq-1} is an approximation of the full system dynamics with its external disturbances in a small neighborhood of the operating point.
It can therefore be used to design feedback controllers, taking into account the full knowledge of the local dynamics of the system.
However, the ultra-local nature of this description makes it necessary to continuously update this model by computing a real-time estimate $F_{\rm{est}}$ of $F$ to account for changes in the operating point as well as in external disturbances. Following the procedure described in detail in \cite{mfc1}, this estimate $F_{\rm{est}}$}
can be computed for $\nu = 1$ as\footnote{See \cite{mfc2} for the case $\nu = 2$.}
\begin{equation}\label{estim}
F_{\rm{est}}(t)  =-\frac{6}{\tau^3}\int_{t-\tau}^t \left\lbrack (\tau -2\sigma)y(\sigma)+ \alpha\sigma(\tau -\sigma)u(\sigma) \right\rbrack d\sigma,
\end{equation}
where $\tau > 0$ may be quite small. This quantity can be computed based on the sole knowledge of the applied input $u$ and the measured output $y$.
In practice, the integral in Eq. (\ref{estim}) is usually replaced by a digital filter.

\subsection{Intelligent proportional controllers}

{For $\nu = 1$, we can define the \emph{intelligent proportional controller} (\emph{iP}) \cite{mfc1}} as
\begin{equation}\label{ip}
    u = -\frac{F_{\rm{est}} - \dot{y}^\ast +K_p e }{\alpha}
\end{equation}
where 
\begin{itemize}
    \item $y^\star$ is the reference trajectory; 
    \item $e = y - y^\star$ is the tracking error; 
    \item $K_p \in \mathbb{R}$ is the tuning proportional gain. 
\end{itemize}
Combining Eqs \eqref{eq-1} and \eqref{ip} yields
$$
\dot{e} + K_p e - F + F_{\rm{est}} = 0
$$
By picking $K_p> 0$, if the estimate $F_{\rm{est}}$ of $F$ is good, i.e., $F - F_{\rm{est}} \approx 0$, then $\lim_{t \to \infty} e(t) \approx 0$. This is equivalent to saying that the output $y$ approximately reaches the reference $y^*$.

\section{A model}\label{descrip}
\subsection{Short description}
Consider (see, e.g., \cite{strugarek2019use})
\begin{equation*}\label{model}
\begin{cases}
\dot x_1=\beta_E x_3 (1-\frac{x_1}{\mathcal{K}})-(\nu_E+\delta_E)x_1\\
\dot x_2=(1-\nu) \nu_E x_1 - \delta_M x_2\\
\dot x_3=\nu\nu_E x_1 \frac{x_2}{x_1+\gamma_S x_4}-\delta_F x_3\\
\dot x_4=u-\delta_S x_4\\
\end{cases}
\end{equation*}
where $\mathcal{K}$ is the carrying capacity for the eggs. The states $x_1$, $x_2$, $x_3$, and $x_4$ correspond, respectively, to $\mathcal{E}$, $\mathcal{M}$, $\mathcal{F}$, and $\mathcal{M}_s$: 

\begin{itemize}
    \item $\mathcal{E}$ is the quantity of viable eggs (in fact, of mosquitoes in the aquatic phase, since we don't add extra compartments for pupae and larvae);
    \item $\mathcal{F}$ is the quantity of females that were fertilized by wild males (and which, therefore, are able to lay viable eggs);
    \item $\mathcal{M}$ is the quantity of wild males;
    \item $\mathcal{M}_s$ is the quantity of sterile males.
\end{itemize}


\subsection{Continuous-time control setting}
Introduce
\begin{equation*}\label{model}
\begin{cases}
\dot x_1=\beta_E x_3 (1-\frac{x_1}{\mathcal{K}})-(\nu_E+\delta_E)x_1\\
\dot x_2=(1-\nu) \nu_E x_1 - \delta_M x_2\\
\dot x_3=\nu\nu_E x_1 \frac{x_2}{x_1+\gamma_S V}-\delta_F x_3\\
\end{cases}
\end{equation*}
$\beta_E$ represents the rate at which females lay eggs; $\mathcal{K}$ represents the carrying capacity of the environment for the mosquitoes in the aquatic phase; $\nu_E$ is the transition rate from the aquatic phase to the adult phase; $\delta_E, \delta_F, \delta_M$ and $\delta_S$ are, respectively, the death rates for the aquatic phase, the (fecundated) females, the wild males, and the sterile males; $\nu$ is the probability that a mosquito is born female, and $\gamma_S$ is a parameter that represents the possible mating preferences of females between wild and sterile males.

The continuous auxiliary control variable $V = x_4$ yields 
the following ultra-local model, where $y = x_1$,
$$\dot y = F+\alpha V$$
The corresponding intelligent proportional controller reads
$$V=\frac{-F_{\rm{est}} +K_p e+\dot y^\star}{\alpha}$$
where, like in Eq. \eqref{ip}, $F_{\rm{est}}$ is an estimate of $F$, $y^\star$ is the reference trajectory, and $e = y-y^\ast$ is the tracking error.
\section{Discrete and realistic control}\label{algo}

In practice, sterile male mosquitoes are not released continuously (all the time) but only at certain moments—typically twice a week for Aedes mosquitoes (see \cite{gato2021sterile}). Therefore, it is natural to work in an impulsive control setting.

We are looking for $u$ as a sequence of pulses every $J$ days such that $\text{mean}|x_4-V|$ is minimal during the period. 
Thus, considering that the sampling period is 1 day, the discrete control $u$ is of the form $\begin{bmatrix}\delta_0&0&...&0& \delta_1& 0&...&0&\delta_k&...\end{bmatrix}$
at times $\begin{bmatrix}1& 2&...&J-1& J &J+1&... &kJ-1& kJ&...\end{bmatrix}$.
Where $\delta_k$ is a pulse of amplitude $\delta_k$.\\
Thanks to the principle of superposition and using the dynamic $\dot x_4=u-\delta_S x_4$, which is assumed to be relatively well known, the amplitude $\delta_k$ is determined as follows 
$$\delta_k=\frac{V(k)-\text{mean}(\mathcal{I}_m({k+1}...{k+J}))}{\text{mean}(\mathcal{I}(1...{J}))}$$
where $V(k)$ is the value of the continuous control $V$ on the $k$th day, $\mathcal{I}_m({k+1}...{k+J})$ is the set of impulse responses of the dynamics corresponding to $x_4$ (transfer function $\frac{1}{s+\delta_S}$ ) at input $u$ for $0<t<k$, padded with zeros for $k\leq t<k+J$. 
$\mathcal{I}(1...{J})$ is the $x_4$'s response to a pulse of amplitude 1.

\section{Results}\label{sec2}

\subsection{Generalities}

As one can expect, for a vector-borne disease, the risk of transmission goes to $0$ when the number of vectors goes to $0$. This is particularly clear in the explicit expression for the basic reproduction number, $R_0$, for the Ross-MacDonald model (see, e.g., \cite{AndersonMay})

\begin{equation*}
    R_0 = \frac{\beta^2 p p' V }{\alpha \mu H} ,
\end{equation*}

\noindent where $\beta$ is the biting rate of the vector, $p$ is the chance that a bite transmits the disease from the vector to the host, $p'$ is the chance that the bite transmits the disease from the host to the vector, $V$ is the vector population, $H$ is the human population, $
\alpha^{-1}$ is the average infectious period for humans, and $
\mu^{-1}$ is the average life expectancy of a vector (in this model, once infectious, vectors are assumed to remain infectious until they die).

Placing ourselves in the setting where the epidemic spread has still not begun, so that the basic reproduction number, $R_0$, is the appropriate quantity (if not, one can choose to use the effective reproduction number). Then, if 

\begin{equation*}
    V \leq V_c := \frac{\alpha \mu H}{\beta^2 p p'} ,
\end{equation*}

\noindent we will have that $R_0 \leq 1$ and that there will be no epidemic.

Therefore, to avoid epidemic risk, we do not need to reduce the population to zero. It suffices that the number of vectors becomes smaller than $V_c$. This will be the objective of our intervention strategy.

Such a strategy will also have to take into account the logistical constraints of the production and release of sterile mosquitoes. Thus, we cannot expect to release an infinite number of sterile mosquitoes, and we even have an interest in releasing a small number of them when the number of wild mosquitoes is small so that we can use the rest of the available sterile mosquitoes to intervene in other places (for instance, in the setting of a rolling carpet strategy \cite{almeidaRolling2023,almeidaRolling2025}). 

These are the considerations that lead us to define a target trajectory for the wild population, which should lead it to a value below $V_c$ in a finite time, thanks to our sterile male releases.


\subsection{Nominal case}
\begin{figure*}[!ht]
\centering
\subfigure[\footnotesize State evolution]
{\epsfig{figure=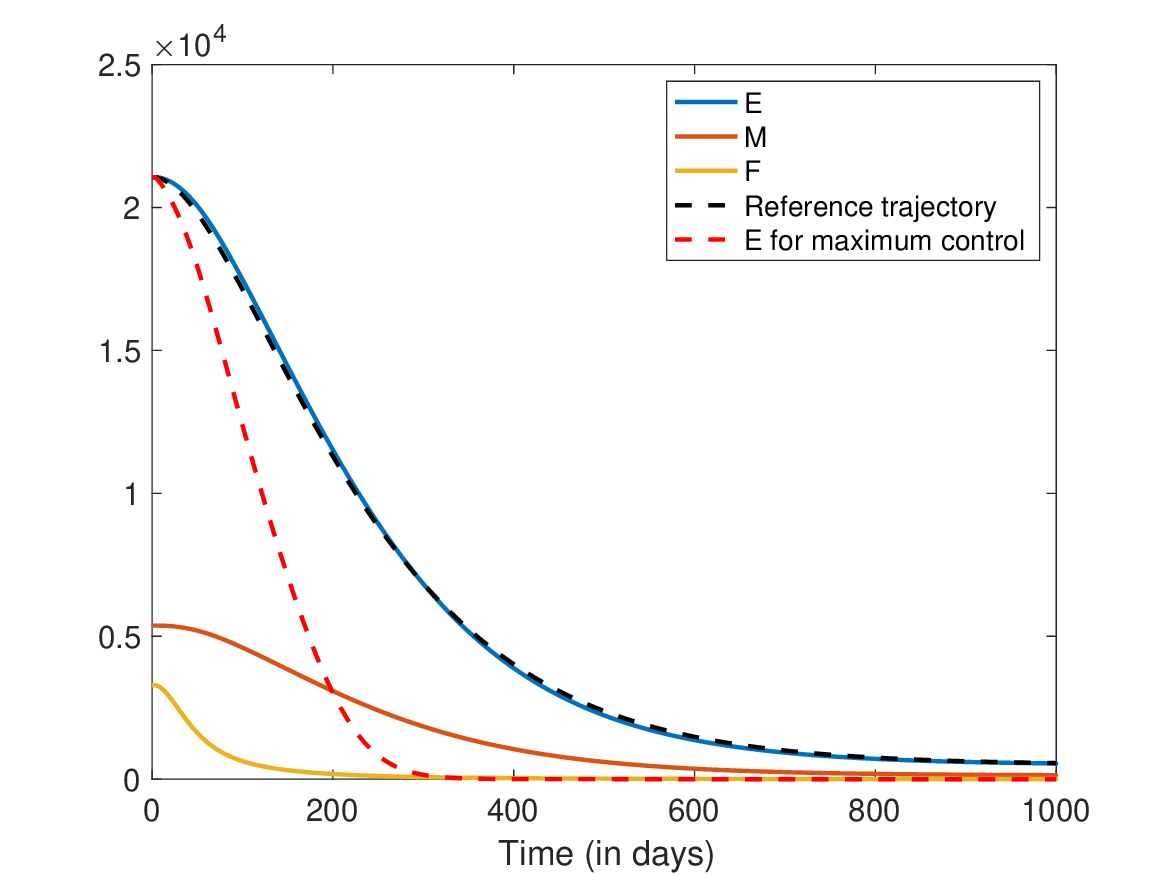,width=0.87\textwidth}}
\\
\subfigure[\footnotesize Continuous control]
{\epsfig{figure=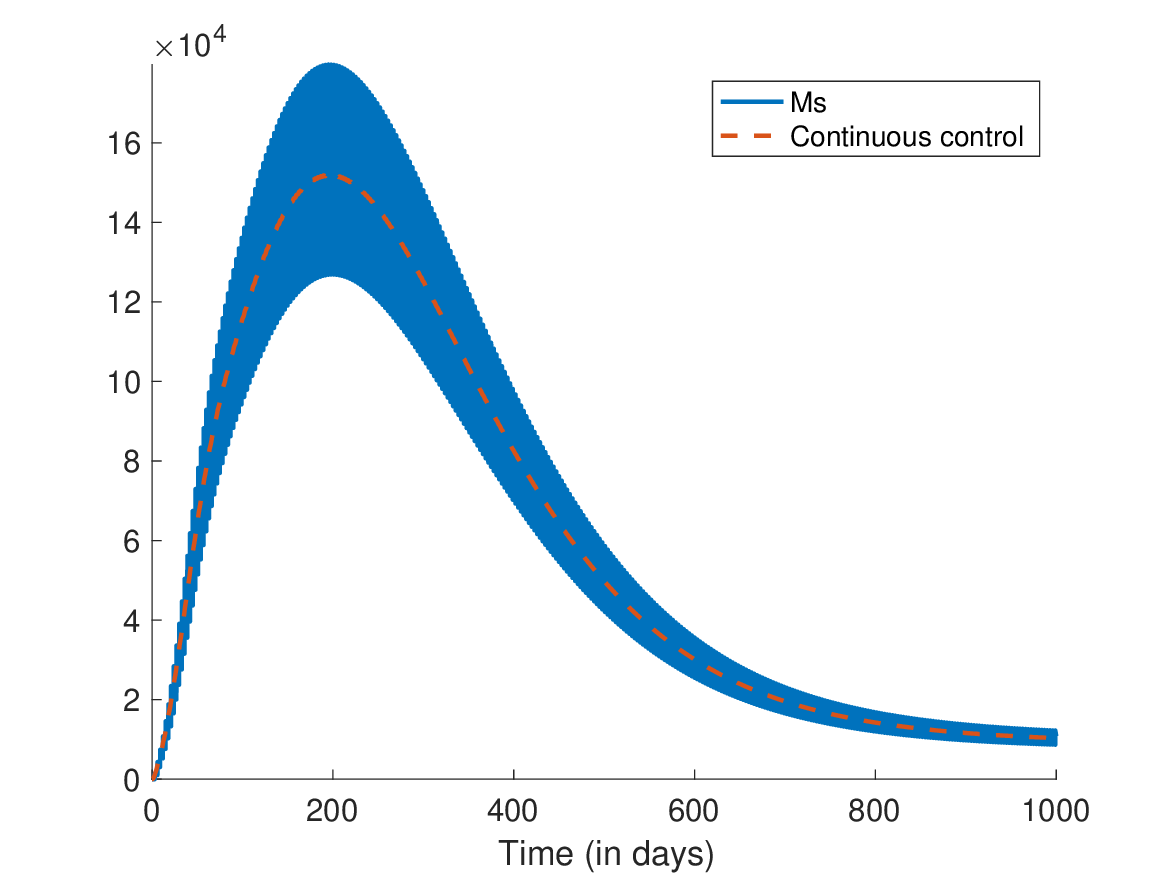,width=0.5\textwidth}}
\subfigure[\footnotesize Impulse control]
{\epsfig{figure=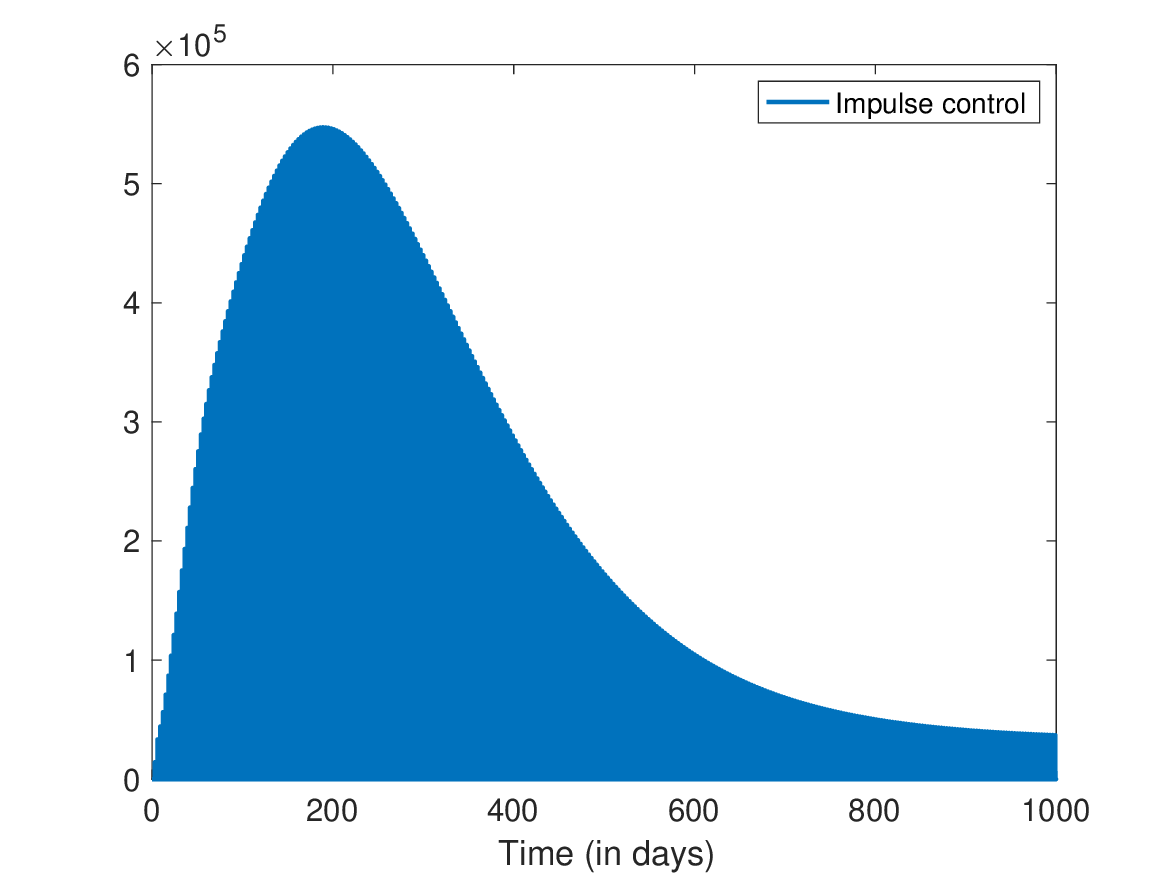,width=0.5\textwidth}}
\caption{Nominal case}\label{S1}
\end{figure*}

Nominal values of considered parameters are 
$\beta_E=10$, $K=22200$, $\nu_E=0.05$, $\delta_E=0.03$, $\nu=0.49$, $\delta_M=0.1$, $\gamma_S=1$, $\delta_F=0.04$, $\delta_S=0.12$ as in \cite{strugarek2019use}.
The initial state value corresponds to the equilibrium point for $V=0$.
%
We consider that pulse control is applied every 3 days, $J=3$ (see Fig. \ref{S1}-(c)). This control is computed using the algorithm explained in Sect. \ref{algo} and, as illustrated in Fig. \ref{S1}-(b), the reconstruction of the continuous control is satisfactory. The reference trajectory is tracked perfectly, as shown in Fig. \ref{S1}-(a). Furthermore, the pulse control is saturated $0<u<10^6$ and if $u=10^6$ at any time, we obtain the dotted curve in Fig. \ref{S1}-(a). Obviously, this last case has an economic cost that is difficult to bear, not to mention the availability of irradiated mosquitoes.

\subsection{Robustness}
\subsubsection{Robustness with respect to $J$}
\begin{figure*}[!ht]
\centering
\subfigure[\footnotesize State evolution]
{\epsfig{figure=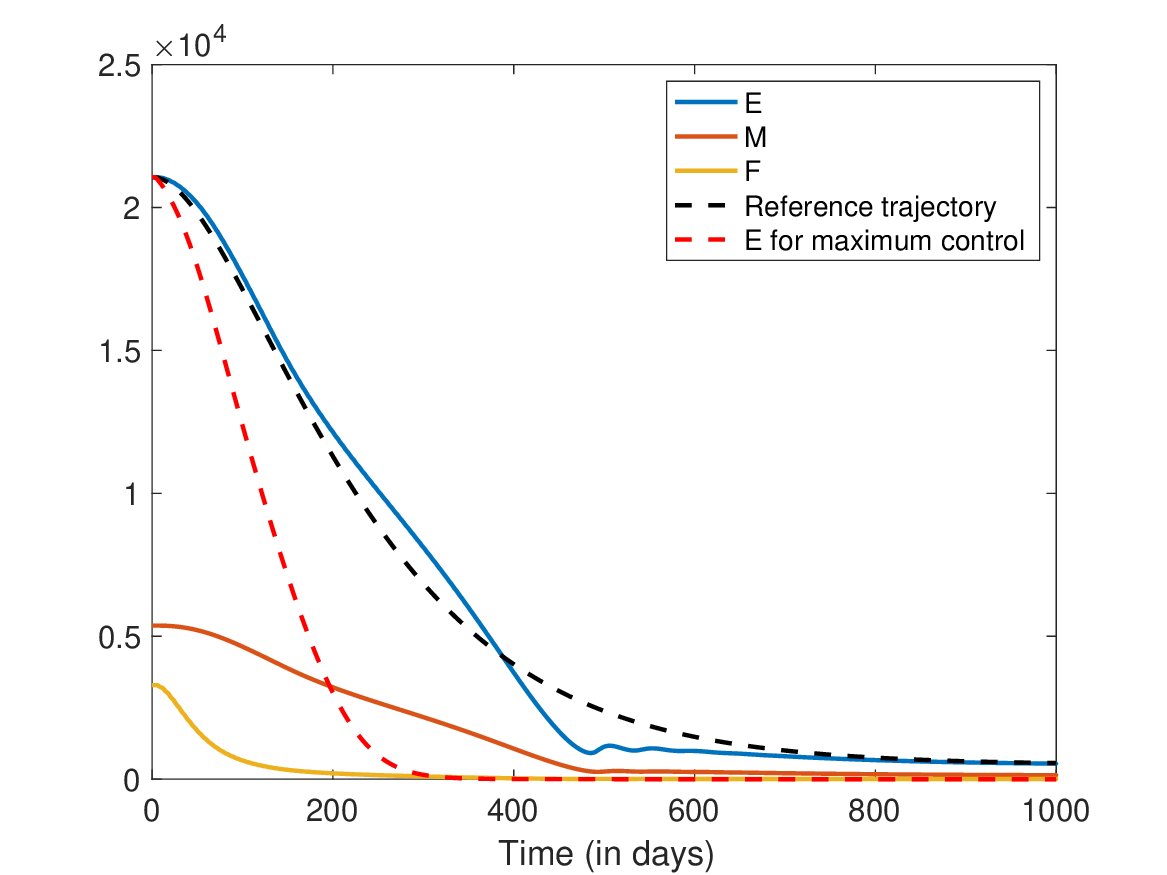,width=0.87\textwidth}}
\\
\subfigure[\footnotesize Continuous control]
{\epsfig{figure=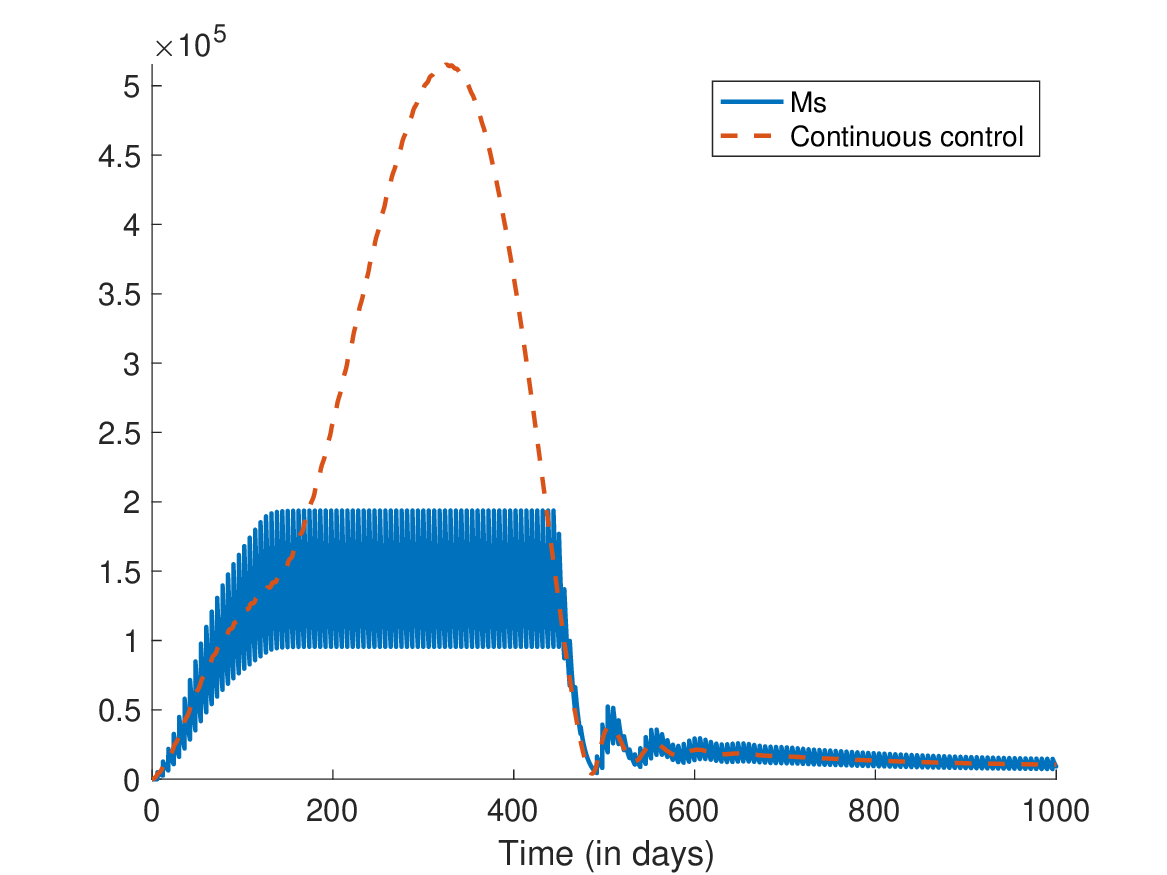,width=0.5\textwidth}}
\subfigure[\footnotesize Impulse control]
{\epsfig{figure=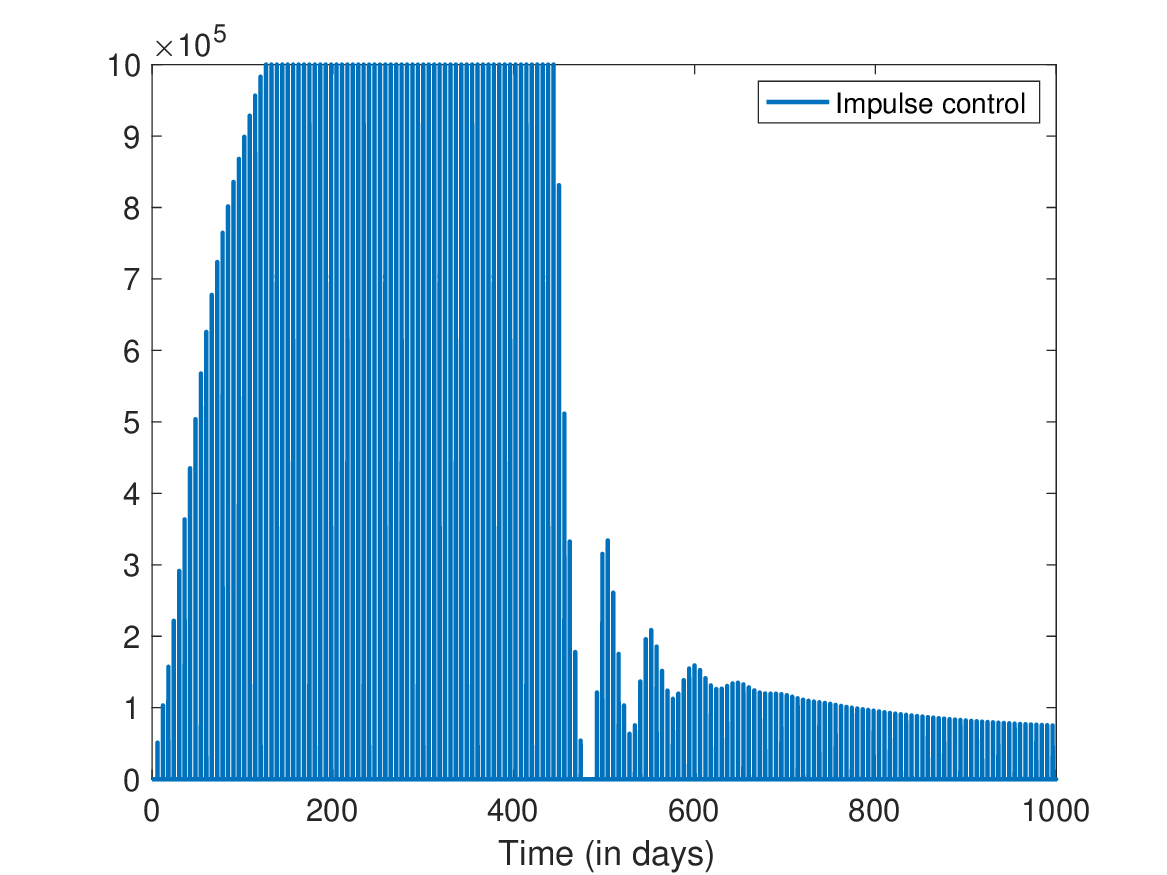,width=0.5\textwidth}}
\caption{Robustness wrt $6$ days instead of $3$}\label{S2}
\end{figure*}

The duration between each treatment or introduction of irradiated mosquitoes, or between each control pulse, has been increased to $J=6$ days. This has degraded the reference trajectory tracking, particularly because saturation is reached in the first period (see Fig. \ref{S2}-(c)). This also explains the deficit in the representation of the continuous command (red dot) by the impulse response of $x_4$ (in blue) in Fig. \ref{S2}-(b). Despite the elements presented above, the trajectory tracking remains very accurate, see Fig. \ref{S2}-(a).

%
\subsubsection{Robustness with respect to the dynamics of $x_4$}
\begin{figure*}[!ht]
\centering
\subfigure[\footnotesize State evolution]
{\epsfig{figure=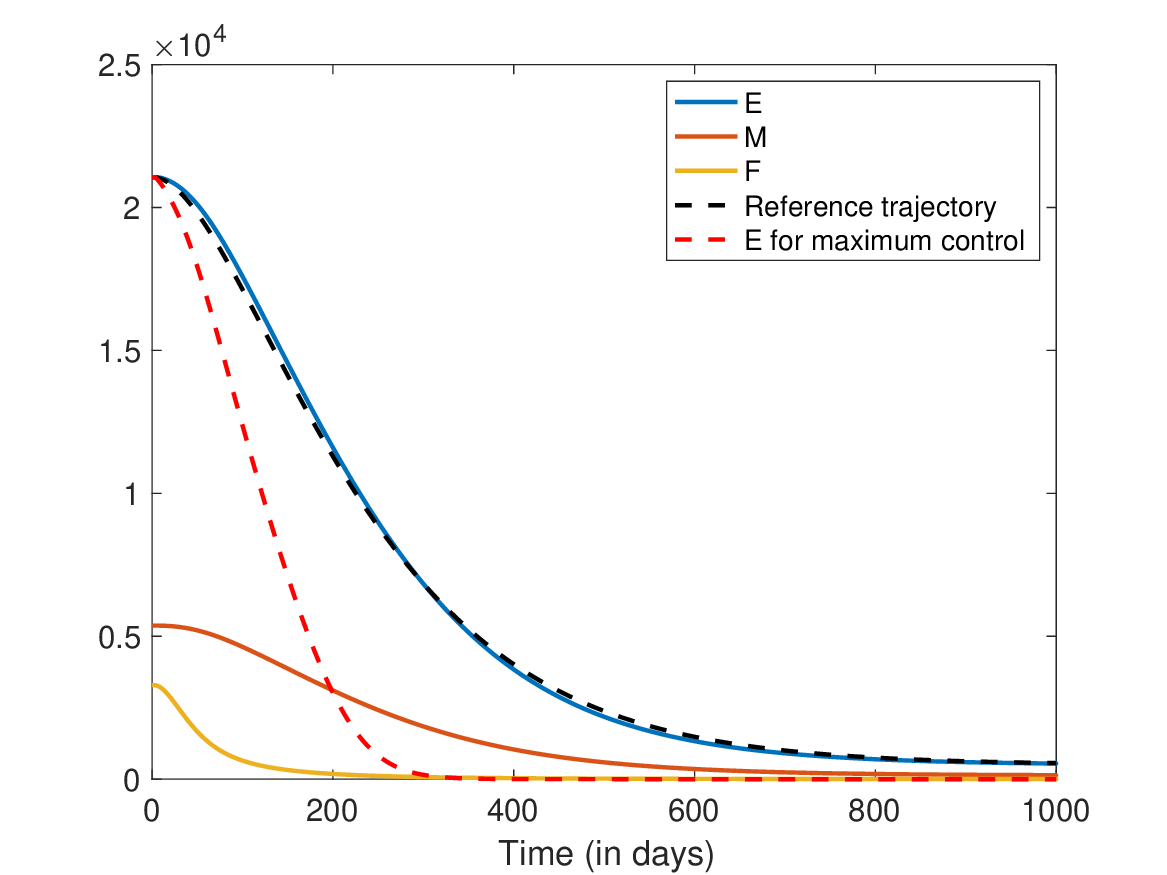,width=0.87\textwidth}}
\\
\subfigure[\footnotesize Continuous control]
{\epsfig{figure=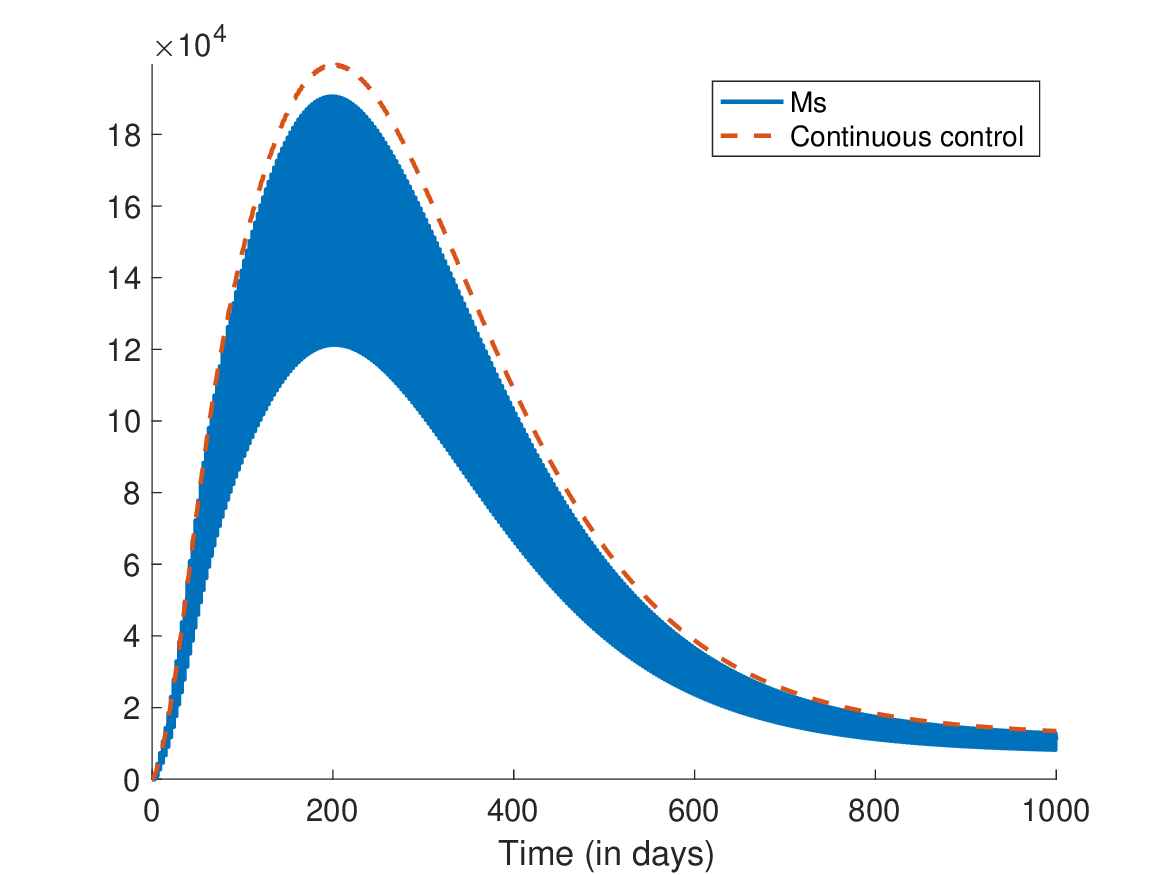,width=0.5\textwidth}}
\subfigure[\footnotesize Impulse control]
{\epsfig{figure=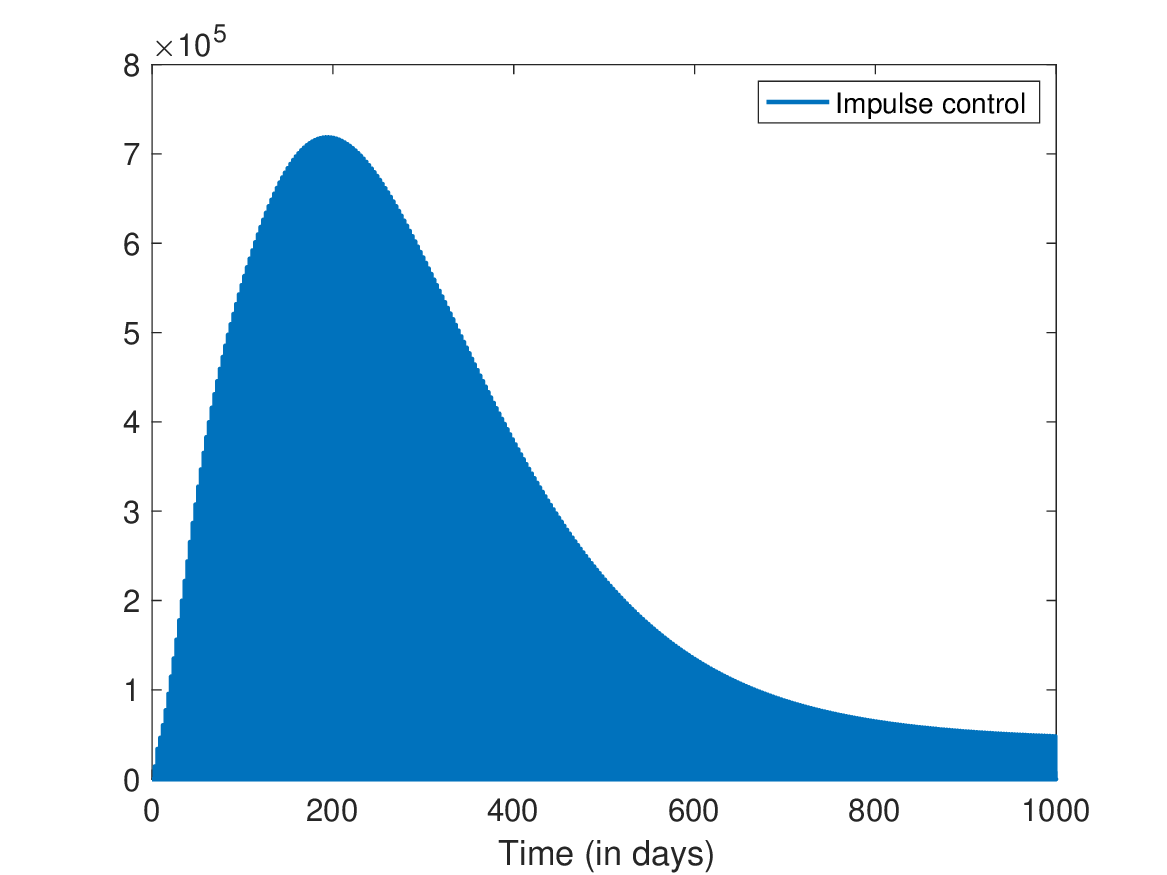,width=0.5\textwidth}}
\caption{Robustness to parametric uncertainties: 100 simulations}\label{S3}
\end{figure*}

Finally, in order to evaluate the robustness of our method, assuming knowledge of the dynamics of $x_4$, we perform a variation of $\delta_S$. 
Thus, to compute the magnitude of control pulse, we consider the nominal value $\delta_S=0.12$ while the model is simulated with $\delta_S*1.3=0.156$. Performance remains very good, see Fig. \ref{S3}-(a)-(c), but obviously, the continuous-discrete representation no longer coincides, as illustrated in Fig. \ref{S3}-(b).

%
\subsubsection{Robustness with respect to model uncertainties}
\begin{figure*}[!ht]
\centering
{\epsfig{figure=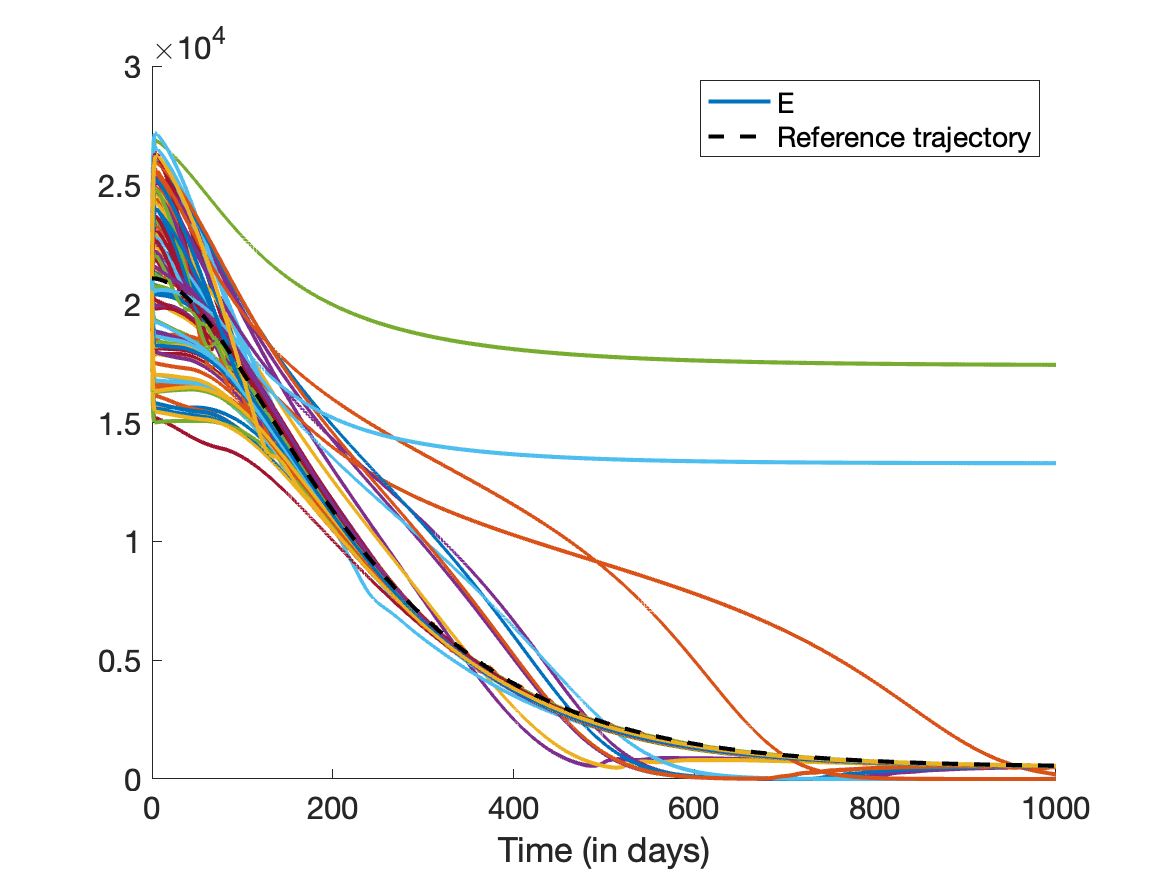,width=0.87\textwidth}}
%
\caption{Robustess to model incertities : 100 simulations}\label{S4}
\end{figure*}

A random draw of $8$ parameters ($\{I_1,...,I_8\}$) between $0.7$ and $1.3$ introduces parametric uncertainty into the model. The nominal parameters $\beta_E=10$, $K=22200$, $\nu_E=0.05$, $\delta_E=0.03$, $\nu=0.49$, $\delta_M=0.1$, $\gamma_S=1$, $\delta_F=0.04$, $\delta_S=0.12$ are therefore multiplied by $I_1$, ..., and $I_8$.\\
\\
$100$ simulations are performed with different random draws. The trajectory tracking is illustrated in Figure \ref{S4}. In only two cases is tracking not guaranteed because control is saturated at all times. This particular combination of random choices is therefore particularly degrading.

\section{Conclusion}\label{conclusion}
Concrete implementations of model-free control with respect to the Sterile Insect Technique, which have often been achieved in diverse domains (see, e.g., \cite{fuel,forest} for recent concrete examples), should not pose any serious challenge. 

We are aware of the impossibility of accurately performing the measurements required for any closed-loop control system. By making our proposal robust against sensor failures such as \cite{longwyIJC}, we could further improve its positioning on the \emph{Technology Readiness Levels} (\emph{TRL}) scale (see, e.g., \cite{trl}), which is well known to be related to system complexity. This is a point that seems quite ignored in today's connections between control and AI.

\end{document}